\documentstyle[twoside,12pt]{article}

\newcommand{\chapter}{\section}

\begin{document}

\newtheorem{Thm}{Theorem}
\newtheorem{Ax}{Axiom}
\newtheorem{Prop}{Proposition}
\newtheorem{Cor}[Prop]{Corollary}
\newtheorem{Main}{}
\renewcommand{\theMain}{}
\newtheorem{Lem}[Prop]{Lemma}
\newtheorem{Fact}{Fact}
\renewcommand{\theFact}{}

\newtheorem{Def}{Definition}
\newtheorem{rmk}{Remark}
\newenvironment{Rmk}{\begin{rmk}\em}{\end{rmk}}
\newtheorem{exm}{Example}
\newenvironment{Exm}{\begin{exm}\em}{\end{exm}}

\newcommand{\qed}{\par {\EM QED} }
\newtheorem{prf}{Proof}
\renewcommand{\theprf}{}
\newenvironment{Prf}{\begin{prf}\em}{\qed\end{prf}}
\newtheorem{prff}{}
\renewcommand{\theprff}{}
\newenvironment{Prff}{\begin{prff}\em}{\qed\end{prff}}





\newcommand{\YES}[1]{#1}
\newcommand{\NOT}[1]{}

\newcommand{\cA}{{\cal A}}
\newcommand{\cB}{{\cal B}}
\newcommand{\cC}{{\cal C}}
\newcommand{\cD}{{\cal D}}
\newcommand{\cE}{{\cal E}}
\newcommand{\cF}{{\cal F}}
\newcommand{\cG}{{\cal G}}
\newcommand{\cH}{{\cal H}}
\newcommand{\cI}{{\cal I}}
\newcommand{\cJ}{{\cal J}}
\newcommand{\cK}{{\cal K}}
\newcommand{\cL}{{\cal L}}
\newcommand{\cM}{{\cal M}}
\newcommand{\cN}{{\cal N}}
\newcommand{\cO}{{\cal O}}
\newcommand{\cP}{{\cal P}}
\newcommand{\cQ}{{\cal Q}}
\newcommand{\cR}{{\cal R}}
\newcommand{\cS}{{\cal S}}
\newcommand{\cT}{{\cal T}}
\newcommand{\cU}{{\cal U}}
\newcommand{\cV}{{\cal V}}
\newcommand{\cW}{{\cal W}}
\newcommand{\cX}{{\cal X}}
\newcommand{\cY}{{\cal Y}}
\newcommand{\cZ}{{\cal Z}}

\newcommand{\bbb}[1]{{\mbox{\bf #1}}}

\newcommand{\bN}{\bbb{N}}
\newcommand{\bZ}{\bbb{Z}}
\newcommand{\bR}{\bbb{R}}
\newcommand{\bC}{\bbb{C}}
\newcommand{\bQ}{\bbb{Q}}
\newcommand{\bT}{\bbb{T}}

\newcommand{\noind}[1]{{\setlength{\parindent}{0cm} #1}}
\newcommand{\parsk}{\par\medskip}

\newcommand{\varend}{\end{document}}

\newcommand{\LP}{\left(}  \newcommand{\RP}{\right)}
\newcommand{\LQ}{\left[}  \newcommand{\RQ}{\right]}
\newcommand{\LB}{\left\{} \newcommand{\RB}{\right\}}
\newcommand{\LA}{\left\langle} \newcommand{\RA}{\right\rangle}
\newcommand{\LS}{\left|}  \newcommand{\RS}{\right|}
\newcommand{\LN}{\left\|} \newcommand{\RN}{\right\|}
\newcommand{\bLP}{\bigl(}  \newcommand{\bRP}{\bigr)}

\newcommand{\ts}{{\otimes}} \newcommand{\TS}{{\bigotimes}}

\newcommand{\es}{\emptyset}
\newcommand{\sm}{\setminus} \newcommand{\st}{\subset}
\newcommand{\eps}{\varepsilon}

\newcommand{\beware}{~$\Psi\!\Psi\!\Psi$~}

\newcommand{\EM}{\em\bf}

\newcommand{\BE}{\begin{equation}}  \newcommand{\EE}{\end{equation}}
\newcommand{\BER}{$$\begin{array}}  \newcommand{\EER}{\end{array}$$}
\newcommand{\head}[1]{
             \par\bigskip\noind{{\large\bf#1}\nopagebreak\par}\medskip}
\newcommand{\chead}[1]{\begin{center}
              \par\bigskip{\large\bf#1}\nopagebreak\par\medskip
              \end{center}}

\newcommand{\fillitem}{\L{\embox{}}\vspace{-.6cm}}

\newcommand{\dfrac}[2]{\displaystyle{\frac{#1}{#2}}}
\newcommand{\toprove}[1]{{\buildrel\mbox{?}\over#1}}
\newcommand{\DEF}{{\buildrel\mbox{\scriptsize\,\,def\,\,}\over=}}
\newcommand{\BAR}[1]{\overline{#1}}
\newcommand{\I}{\infty}
\newcommand{\al}{\alpha} \newcommand{\be}{\beta}
\newcommand{\OM}{\Omega} \newcommand{\om}{\omega}
\newcommand{\la}{\lambda}
\newcommand{\mod}{\mbox{ mod }}
\newcommand{\half}{\frac{1}{2}}
\newcommand{\NI}{{n\to\infty}}
\newcommand{\ang}{{<\!\!\!)}}
\newcommand{\arc}[1]{\breve{#1}}

\newenvironment{txeqn}[1]{\begin{equation}\label{#1}\begin{array}{l}}{
                         \end{array}\end{equation}}



\newcommand{\entropy}{{\mbox{entropy}\,}}
\newcommand{\tr}{{\mbox{tr}\,}}
\newcommand{\OR}{{\mbox{ {\EM or}}\,}}

\title{A Note About Entropy}
\author{Eliahu Levy\\
Department of Mathematics\\
Technion -- Israel Institute of Technology,
Haifa 32000, Israel\\
email: eliahu@techunix.technion.ac.il}


\date{}


\maketitle
\begin{abstract}
A  mathematical interpretation\NOT{(formula $(2)$)} of the usual
definition of entropy is given. This formulation makes some
properties of entropy immediate.
\end{abstract}


It seems illuminating to work simultaneously with two
classical scenarios where entropy occurs: the set $S$ of
nonnegative real sequences with sum $1$ (discrete probability
distributions) and the set $O$ of trace-class positive-definite
(Hermitian) operators with trace $1$ on a separable Hilbert space
$\cH$  (the quantum-mechanical counterpart of the former). In both
cases the discussed set is a subset of a ``$\|\;\|_1$'' Banach
space (the sequence space $\ell_1$ for $S$ and the space of
trace-class operators for $O$) with its norm $\|\;\|_1$, and is a
subset of a ``$\|\;\|_\I$'' Banach space (the sequence space
$\ell_\I$ for $S$ and the space $B(\cH)$ for $O$) with its norm
$\|\;\|_\I$.
\parsk

In both cases the entropy of an $\om\in S \OR O$ is defined as
usual. For $S$
$$\entropy(\om)=-\sum_i\om_i\ln(\om_i),$$
and for $O$ one applies the same formula for the sequence
$(\om_i)$ of eigenvalues.
\parsk

An ``interpretation'' for this expression may be given as follows:
\parsk

For $\om\in S\OR O$, Define $U(\om)\subset\bR^+\times\bR^+$
(where $\bR^+=[0,\I[\;\subset\bR$) as follows: a pair
$(r_1,r_\I)$, $r_1,r_\I\ge0$, belongs to $U(\om)$ iff
for two members of the $\|\;\|_1$-Banach space with difference
$\om$, the open $\|\;\|_1$-ball of radius $r_1$ centered in one
of them and the open $\|\;\|_\I$-ball of radius $r_\I$
centered in the other {\em do not intersect}.
\parsk

Clearly,
$$(r_1',r_\I')\le(r_1.r_\I)\in
U(\om)\Rightarrow(r_1',r_\I')\in U(\om).$$

In other words, the boundary $\partial U(\om)$ in
$\bR^+\times\bR^+$ is the ``graph'' a of a non-increasing function
$F:r_1\leftrightarrow r_\I$ (where in case of a jump the
vertical segment spanning the jump is to be added to this
``graph'') and $U(\om)$ is the set of points in the first
quadrant below this ``graph''.
\parsk

In fact, it is easy to find $U(\om)$ explicitly: one easily
checks for $S$, and, somewhat surprisingly, equally easily for
$O$, that in both cases $F$ is given by:
$$r_1=\sum_i\max(0,\om_i-r_\I),\leqno(1)$$
where in the $O$ case the $\om_i$ are the eigenvalues. (In other
words, a pair $(r_1,r_\I)$ belongs to $U(\om)$ iff $r_1$ is
$0$ or less than the expression $(1)$.)
\parsk

Thus, the function $F$ is continuously decreasing from
$r_1=0,r_\I=\max(\om_i)=\|\om\|_\I$ to
$r_\I=0,r_1=1$ and is piecewise linear, namely linear in
intervals between the points where $r_\I$ attains as values
the $\om_i$'s, the slope of $r_\I$ w.r.t.\ $r_1$ in such an
interval being the reciprocal of the number $n$ of the
$\om_i$'s greater than the range of $r_\I$ in this
interval.
\parsk

Now integrate $\ln(r_\I)\,dr_1$ over the whole $\partial
U(\om)$. An interval $[r_1',r_1'']$ spoken about in the last
paragraph, the values of $r_\I$ at its endpoints being
$r_\I''=\om_{i+1}<r_\I'=\om_i$ (we ordered the
$\om_i$ in non-increasing order), will contribute (note that
here $dr_1=n\cdot\,dr_\I$):
$$n[\om_{i+1}\ln(\om_{i+1})-\om_i\ln(\om_i)]+r_1''-r_1',$$
and adding up all the intervals one get the sought-for expression for
the entropy:
$$\entropy(\om)=-\int_{\partial
U(\om)}\ln(r_\I)\,dr_1-1.\leqno(2)$$

An immediate consequence of $(2)$ is that any linear
transformation of sequences (case $S$) or operators (case $O$),
that both does not increase $\|\;\|_\I$ and does not increase
$\|\;\|_1$, is entropy non-decreasing (it enlarges balls, thus
shrinks $U(\om)$, hence lowers $\partial U(\om)$). In
particular, one deduces the fact, important in quantum theory,
that for any fixed basis in the Hilbert space $\cH$, the operation
of erasing all off-diagonal entries in the matrix of an $\om\in O$
is entropy non-decreasing. More generally, let $X$ be a closed real
linear subspace of the real (or Hermitian) elements of the
$\|\;\|_1$-Banach space (either for $S$ or for $O$), so that
$S_X=X\cap S\OR O_X=X\cap O$ spans $X$. Suppose any $\om\in S\OR O$ has
a ``conditional expectation'' (necessarily unique) $\om_X\in S_X\OR O_X$,
in the sense that
$$\forall a\in X\,\,\LA\om_X,a\RA=\LA\om,a\RA,$$
where in the $S$-case $\LA a,b\RA:=\sum a_ib_i$ and in the $O$-case
$\LA a,b\RA:=\tr(ab)$. Then one has direct-sum decompositions of the
real $\|\;\|_1$- and $\|\;\|_\I$-Banach spaces into the closures of $X$ and
of its annihilator, and the projection $\pi_X$ on $X$ in these
decompositions is norm-non-increasing for both norms. (Indeed, write each
real/Hermitian $a$ in the $\|\;\|_1$-space as $\al_+\om_+-\al_-\om_-$,\,\,
$\al_+,\al_-\ge0$,\,\,$\al_++\al_-=\|a\|_1$,\,\,$\om_+,\om_-\in S\OR O$.
Then $\|\pi_X(a)\|_1=\|\al_+{\om_+}_X-\al_-{\om_-}_X\|_1\le\al_++\al_-=\|a\|_1$.
For another real/Hermitian $b$, $\LA\pi_X(a),b\RA=\LA a,\pi_X(b)\RA$ hence
$|\LA\pi_X(a),b\RA|\le\|a\|_\I\|\pi_X(b)\|_1\le\|a\|_\I\|b\|_1$ implying
$\|\pi_X(a)\|_\I\le\|a\|_\I$.)
And one deduces that taking the above conditional expectation is entropy
non-decreasing.

\end{document}